\documentclass[12pt]{amsart}
\usepackage{amsfonts,amsmath,amssymb,amsthm,booktabs,float,geometry,enumerate}
\newtheorem{theorem}{Theorem}[section]
\newtheorem{prop}[theorem]{Proposition}
\newtheorem{lem}[theorem]{Lemma}
\newtheorem{cor}[theorem]{Corollary}
\theoremstyle{definition}
\newtheorem{rem}[theorem]{Remark}
\newtheorem{mydef}[theorem]{Definition}
\newtheorem{example}[theorem]{Example}

\begin{document}
\title{Poincar\'e polynomials of nilmanifolds associated with path graphs}
\author[Aldi]{Marco Aldi}
\author[Bevins]{Samuel Bevins}
\author[Chatzitheodoridis]{Eleftherios Chatzitheodoridis}
\author[Da Silva]{Sergio Da Silva}
\author[Frias]{Quincy Frias}
\author[Mej{\'i}a G{\'o}mez]{Tom\'as Mej{\'i}a G{\'o}mez}
\subjclass{57R19 (Primary), 05C38, 17B55, 17B56 (Secondary).}
\keywords{Dani--Mainkar nilmanifolds, path graphs, Betti numbers, cohomology}
\address[Aldi]{Department of Mathematics and Statistics, Virginia Commonwealth University, Grace E. Harris Hall, 1015 Floyd Ave., Richmond, VA 23284, USA}
\email{maldi2@vcu.edu}
\address[Bevins]{Department of Physics, William \& Mary, Small Hall, 300 Ukrop Way, Williamsburg, VA 23185, USA}
\email{sjbevins@wm.edu}
\address[Chatzitheodoridis]{Department of Mathematics and Statistics, Virginia Commonwealth University, Grace E. Harris Hall, 1015 Floyd Ave., Richmond, VA 23284, USA}
\email{chatzitheode@vcu.edu}
\address[Da Silva]{
    Department of Mathematics and Economics, 
    Virginia State University, 
    1 Hayden Dr.,
    Petersburg, VA 23806, USA}
\email{sdasilva@vsu.edu, smd322@cornell.edu}
\address[Frias]{Department of Mathematical Sciences, George Mason University, Exploratory Hall, 4400 University Dr., Fairfax, VA 22030, USA}
\email{jfrias3@gmu.edu}
\address[Mej{\'i}a G{\'o}mez]{
    Fakultät für Mathematik, Ruhr-Universität Bochum, Universitätsstraße 150, 44780 Bochum, Germany
}
\email{Tomas.MejiaGomez@rub.de}

\begin{abstract}
We give explicit and recursive formulas for the Poincar\'e polynomials of Dani--Mainkar nilmanifolds associated with path graphs. We also prove that these polynomials are unimodal, and we calculate a range of Betti numbers asymptotically. As an application, we derive a sharp lower bound for the dimension of the total de Rham cohomology of a Dani--Mainkar nilmanifold.
\end{abstract}
\maketitle

\section{Introduction}\label{sec:1}

In their seminal work \cite{DaniMainkar05}, Dani and Mainkar associate with any finite simple graph $G$ a 2-step nilmanifold $\mathcal N(G)$. In addition to the original motivation of producing examples of nilmanifolds with special properties, this construction has been generalized to hypergraphs by the first and second author \cite{AldiBevins24}, and it has been used to prove results about general 2-step Lie algebras and hypergraphs by the first author, Gabrielsen, Grandini, Harris, and Kelley \cite{AGGHK, AGGHK2}.

From a topological perspective, it is natural to investigate the de Rham cohomology of $\mathcal N(G)$. By a theorem of Nomizu \cite{Nomizu54}, the calculation of the de Rham cohomology of $\mathcal N(G)$ is in principle equivalent to a linear algebra problem, namely calculating the kernel of a finite-dimensional operator, the Lie algebra Laplacian. However, since the dimension of the underlying vector space grows exponentially with the dimension of $\mathcal N(G)$, this approach does not make the problem more tractable. Nevertheless, de Rham cohomology is understood in small degrees for all graphs in work of Pouseele and Tirao \cite{PouseeleTirao09} and of the first author, Butler, Gardiner, Grandini, Lichtenwalner, and Pan \cite{ABGGLP}. It is also understood in all degrees for specific families of graphs, including complete graphs in the work of Kostant \cite{Kostant61} and Sigg \cite{Sigg96} and star graphs in the work of Armstrong, Cairns, and Jessup \cite{ArmstrongCairnsJessup97} and of the first author, Butler, Gardiner, Grandini, Lichtenwalner, and Pan \cite{ABGGLP}.

In this paper, we extend this line of work by calculating the Betti numbers of the 2-step nilmanifold $\mathcal N(P_n)$ for all $n$, where $P_n$ denotes the path graph with $n$ vertices. It was shown by the first author, Butler, Gardiner, Grandini, Lichtenwalner, and Pan \cite{ABGGLP} that the Betti numbers of $\mathcal N(G)$ can be thought of as a weighted count of the induced subgraphs of $G$ with weights given by the so-called \emph{essential Betti numbers} of $G$. This observation is particularly useful when the structure of the induced subgraphs is relatively tame. Path graphs provide a testing ground for this approach, for their induced subgraphs are all disjoint unions of smaller paths.

Following this strategy, our key technical result is a recursive formula for the essential Betti numbers of path graphs. This goal is achieved through a careful analysis of the Dixmier long exact sequence obtained by removing a vertex at one end of the path, as developed by Dixmier \cite{Dixmier55}, Pouseele and Tirao \cite{PouseeleTirao09}, and the first and second author \cite{AldiBevins24}.

Out of the resulting recursive formulas for the essential Betti numbers of path graphs, we draw a recursion for the (essential) Poincar\'e polynomials of $\mathcal N(P_n)$, which we solve to extract an explicit formula. We also show that the Betti numbers of $\mathcal N(P_n)$ are unimodal, a property that is established for the Betti numbers of various classes of nilmanifolds by Armstrong \cite{Armstrong97}.

We then use our explicit formulas to derive asymptotic information regarding the de Rham cohomology of $\mathcal N(P_n)$. Specifically, we are able to show that the dimension of the total cohomology of $\mathcal N(P_n)$ grows exponentially with $n$. As a result, we obtain a lower bound for the dimension of the total cohomology of $\mathcal N(G)$ in terms of the diameter of $G$. We also describe the asymptotics of the $(n+k)$-th Betti number of $\mathcal N(P_n)$ for $k \geq 0$.

\subsection*{Organization of the paper} Section \ref{sec:2} comprises the background that we need in this paper. In Section \ref{sec:3}, we produce a recursive formula for the essential Betti numbers of path graphs. In Section \ref{sec:4}, we develop and explicitly solve a recursion for the essential Poincar{\'e} polynomials of path graphs. We perform the same task for the Poincar{\'e} polynomials of path graphs in Section \ref{sec:5}. Lastly, in Section \ref{sec:6}, we study the asymptotic behavior of the total cohomology, as well as a range of Betti numbers for path graphs.

\subsection*{Acknowledgments}
This work was supported by the 4-VA Collaborative Research Grant ``Collaborative Workshops in Topology.'' We thank Rebecca Field and J.D. Quigley for organizing the collaborative workshop in the summer of 2025 during which the bulk of our results was developed. The fourth author's research is supported in part by NSF LEAPS-MPS Grant 2532757. Claude Opus 5.5 and ChatGPT 6 Pro were used to obtain feedback regarding the exposition of preliminary drafts of this manuscript. All the mathematical ideas and the original content are solely due to the authors.

\section{Background}\label{sec:2}
Firstly, we gather the preliminary results that we rely on throughout this paper.

\subsection{The Dani--Mainkar construction}\label{sec:2.1} We recall the Dani--Mainkar construction for a finite simple graph.

\begin{mydef}[{\cite[\S 1]{DaniMainkar05}}]\label{def:1}
Let $G$ be a finite simple graph with vertices $V(G)=\{1,\ldots,n\}$ and edges $E(G)\subseteq \binom{V(G)}{2}$. The \emph{Dani--Mainkar Lie algebra} of $G$ is the 2-step nilpotent Lie algebra $\mathcal L(G)$ with underlying real vector space spanned by $x_1,\ldots,x_n$ and $x_{i,j}$ for $1\le i<j\le n$ whenever $\{i,j\}\in E(G)$. The Lie bracket is uniquely defined by setting \[[x_i,x_j]=\begin{cases}
x_{i,j}, & \textrm{if }\quad \{i,j\}\in E(G);\\
\\
0, & \textrm{otherwise}.
    \end{cases}\]
\end{mydef}

\begin{example}\label{ex:2}
Let $P_n$ be the path graph on $n$ vertices. Then, $\mathcal L(P_n)$ is the real vector space of dimension $2n-1$ generated by $x_1,\ldots,x_n,x_{1,2},x_{2,3},\ldots,x_{n-1,n}$, with Lie bracket such that $[x_i,x_{i+1}]=x_{i,i+1}$ for all $i=1,\ldots,n-1$. 
\end{example}

Let $\mathcal N(G)$ be a nilmanifold obtained as the quotient of the Lie group $\exp(\mathcal L(G))$ by a discrete and co-compact subgroup $\Gamma$. By a theorem of Nomizu \cite[Theorem 1]{Nomizu54}, the de Rham cohomology of $\mathcal N(G)$ is isomorphic to the cohomology of the Chevalley-Eilenberg complex $\mathcal C^\bullet(G)$ of $\mathcal L(G)$. It is convenient to realize $\mathcal C^\bullet(G)$, which is defined as the exterior algebra of the dual of $\mathcal L(G)$, as the space of skew-symmetric polynomials in the dual variables $x_i^*,x_{j,k}^*$, with $i,j,k\in\{1,\ldots,n\}$, $j<k$, and $\{j,k\}\in E(G)$. In this notation, the Chevalley-Eilenberg differential is identified with the differential operator 
\[Q=\sum_{\{i,j\}\in E(G), \, i<j}x_i^*x_j^*\frac{\partial}{\partial x_{i,j}^*}.\]

\subsection{Essential cohomology}\label{sec:2.2} The essential cohomology of Dani--Mainkar nilmanifolds was defined by the first author, Butler, Gardiner, Grandini, Lichtenwalner, and Pan.

\begin{mydef}[{\cite[Definition 15]{ABGGLP}}]\label{def:3} Let $G$ be a finite simple graph, and let $S$ be a subset of the vertices of $G$. We denote by $\mathcal C^\bullet(G,S)$ the subcomplex of $\mathcal C^\bullet(G)$ spanned by monomials of the form $x_{i_1}^*\cdots x_{i_p}^*x_{r_1,s_1}^*\cdots x_{r_q,s_q}^*$, where $i_1,\ldots,i_p,r_1,s_1,\ldots,r_q,s_q$ account for all of the elements of $S$, possibly with repetition.

\begin{enumerate}[(I)]
    \item The \emph{$k$-th essential cohomology} $\mathcal H^k(G)$ is the $k$-th cohomology of $\mathcal C^\bullet(G,V(G))$.
    \item The \emph{$k$-th essential Betti number} $\beta_k(G)$ is the dimension of $\mathcal H^k(G)$.
\end{enumerate}

\end{mydef}

Since $\mathcal C^\bullet(G)$ decomposes as the direct sum of the subspaces $\mathcal C^\bullet(G,S)$ as $S$ ranges over all subsets of $V(G)$, we obtain that
\begin{equation}\label{eq:3}
b_k(G)=\sum_{|V(H)|\le 2k-1} \left|\binom{G}{H}\right|\beta_k(H),
\end{equation}
where the sum is over all isomorphism types of graphs $H$ with at most $2k-1$ vertices that occur as induced subgraphs of $G$, and $\binom{G}{H}$ denotes the set of all induced subgraphs of $G$ that are isomorphic to $H$. Furthermore, if we write $+$ for the disjoint union of graphs, then we obtain that $\mathcal L(G_1+G_2)=\mathcal L(G_1)\oplus \mathcal L(G_2)$, which yields
\begin{equation}\label{eq:4}
\beta_k(G_1+G_2)=\sum_{i=0}^k\beta_{i}(G_1)\beta_{k-i}(G_2).
\end{equation}
We denote the Poincar\'e polynomial of $\mathcal N(G)$ by
\[
b(G,x)=\sum_{k=0}^\infty b_k(G) x^k.
\]
We define the \emph{essential Poincar\'e polynomial} of $\mathcal N(G)$ by
\[
\beta(G,x)=\sum_{k=0}^\infty \beta_k(G)x^k,
\]
and we obtain that
\begin{equation}\label{eq:5}
b(G,x)=\sum_H \left|\binom{G}{H}\right|\beta(H,x),
\end{equation}
where the sum is over all isomorphism classes of graphs $H$. In particular, \eqref{eq:4} yields
\[\beta(G_1+G_2,x)=\beta(G_1,x)\beta(G_2,x).\]

\begin{rem}\label{rem:4}
Since the Euler characteristic of a torus bundle vanishes, we find that $b(G,-1)=0$ for every nonempty finite simple graph $G$. Substituting into \eqref{eq:5} and using induction on the number of vertices, we obtain that $\beta(G,-1)=(-1)^{|V(G)|}$.
\end{rem}

\begin{example}\label{ex:5} 
The Dani--Mainkar nilmanifold $\mathcal N(P_1)$ is a circle. Thus, $H^0(\mathcal N(P_1))$ is spanned by $1$, and $H^1(\mathcal N(P_1))=\mathcal H^1(P_1)$ is spanned by $x_1^*$. As a result, $b(P_1,x)=1+x$, while $\beta(P_1,x)=x$.
\end{example}

\begin{example}\label{ex:6} 
The Dani--Mainkar nilmanifold $\mathcal N(P_2)$ is the \emph{$3$-dimensional Heisenberg nilmanifold}, which can be viewed either as a $T^2$-bundle over the circle $S^{1}$ or as principal $S^1$-bundle over $T^2$. It follows that $H^\bullet(\mathcal N(P_2))$ is spanned by $1$, $x_1^*$, $x_2^*$, $x_1^*x_{1,2}^*$, $x_2^*x_{1,2}^*$, and $x_1^*x_2^*x_{1,2}^*$. We infer that $\beta(P_2,x)=2x^2+x^3$, while $b(P_2,x)=1+2x+2x^2+x^3$.
\end{example}

\begin{rem}\label{rem:8}
Let $H$ be an induced subgraph of $G$. Since every induced subgraph $H'$ of $H$ is also an induced subgraph of $G$, it follows from \eqref{eq:3} that $b_k(G)\ge b_k(H)$ for every $k$. Equivalently, $b(G,x)-b(H,x)$ is a polynomial with non-negative coefficients. In particular, looking at the dimension of the total cohomology, we obtain the lower bound $b(G,1)\ge b(P_{{\rm diam}(G)+1},1)$, where ${\rm diam}(G)$ denotes the diameter of the graph $G$.
\end{rem}

\subsection{Cohomology operators}\label{sec:2.3} We introduce the shorthand notation $H_n^k=H^k(P_n)$ and $H_n^\bullet=H^\bullet(P_n)$. For $n>1$, let $A_n,B_n,I_n \colon H^\bullet_{n-1}\to H^\bullet_n$ be the linear maps defined as follows. Let $I_n$ be the inclusion map, let $A_n([\omega])=[x_n^*\omega]$, and let $B_n([\omega])=[x_n^*x_{n-1,n}^*\omega]$. Let $\kappa_n$ be the kernel of the multiplication map $x_n^*\cdot \colon H^\bullet_n\to H^\bullet_n$ defined by $x_n^*\cdot [\omega]= [x_n^*\omega]$ for all $\omega\in H^\bullet_n$, and let $\kappa_n^k=\kappa_n\cap H_n^k$. For notational convenience, we formally define $H^\bullet(P_0)$ to be the $1$-dimensional vector space spanned by the generator $1$, and we define $A_1,I_1 \colon H^\bullet_0\to H^\bullet_1$ to be linear such that $I_1([1])=[1]$ is the unit of the algebra $H^\bullet(P_1)$ and $A_1(1)=[x_1^*]$. For later use, we also introduce the auxiliary operator $D_n \colon H^\bullet_{n-1}\to H^\bullet_n$ such that $D_n([\omega])=[x_{n-1,n}^*x_{n-1}^*\omega]$ for all $[\omega]\in H^\bullet_{n-1}$.

\begin{mydef}\label{def:10}
For each $n$, fix a basis $\{[\omega^i_{n-1}]\}$ of $\kappa_{n-1}$ with specified representative polynomials $\omega^i_{n-1}$ and polynomials $\eta^i_{n-1}\in \mathcal C^\bullet(P_{n-1})$ such that $x_{n-1}^*\omega_{n-1}^i=Q(-\eta_{n-1}^i)$ for each $i$. Let $C_n \colon \kappa_{n-1}\to H^\bullet_n$ be the linear map such that
\[
C_n([\omega^i_{n-1}])=[x_n^*\eta^i_{n-1}+x_{n-1,n}^* \omega^i_{n-1}].
\]
\end{mydef}

\begin{rem}\label{rem:12}
While the definition of $C_n$ does not depend on the particular basis of $\kappa_{n-1}$, it depends on the polynomials representing the cohomology classes of the basis. Suppose that $\omega^i_{n-1}+Q(\xi)$ is used to represent $[\omega_{n-1}^i]$ and that, correspondingly, $\eta_{n-1}^i$ is shifted to $\eta_{n-1}^i-\varphi$ subject to $Q(x_{n-1}^*\xi)=Q(-\varphi)$. We obtain a modified map 
\begin{align*}
C_n([\omega^i_{n-1}+Q(\xi)])& =[x_n^*\eta^i_{n-1}+x_{n-1,n}^* \omega^i_{n-1}]-[x_n^*\varphi-x_{n-1,n}^*Q(\xi)]\\
&=[x_n^*\eta^i_{n-1}+x_{n-1,n}^* \omega^i_{n-1}]-x_n^*[\varphi+x_{n-1}^*\xi].
\end{align*}
Hence, the choices involved in the definition of $C_n$ disappear after quotienting $H^\bullet(P_n)$ by the image of $A_n$. We also see that, for any $[\omega]\in \kappa_{n-1}$, there is a polynomial $\eta$ such that $C_n([\omega])=[x_n^*\eta+x_{n-1,n}^* \omega]$.
In particular, $C_n([\omega])$ can be thought of as a representative of the triple Massey product $\langle x_n^*,x_{n-1}^*,\omega\rangle$.
\end{rem}

\begin{rem}\label{rem:12bis}
For future reference, we can refine Definition \ref{def:10} by requiring that, if $\omega_{n-1}^i\in \mathcal C^\bullet(P_{n-1},S)$, then $\eta_{n-1}^i\in \mathcal C^\bullet(P_{n-1},S\cup\{n-1\})$. We apply this refinement in the context of essential cohomology. Let $\mathcal H_n^k$ be the $k$-th essential cohomology of $P_n$,
let $\mathcal K_n^k=\mathcal H_n^k\cap \kappa_n^k$, and let  $\widetilde{ \mathcal K}_n^k=\mathcal H_n^k\cap \kappa_{n+1}^k$. Then, $C_n (\mathcal K_{n-1}^{k-1}\oplus I_{n-1}( \widetilde{\mathcal K}_{n-2}^{k-1}))\subseteq \mathcal H_n^k$.
\end{rem}

\section{A recursive formula for the essential Betti numbers of $\mathcal N(P_n)$}\label{sec:3}
In this section, we develop a sequence of lemmas that leads to a recursive formula for the essential Betti numbers of path graphs. This recursive formula is the foundation of our work in the next sections, and it yields an efficient tabulation of the essential Betti numbers of path graphs.

Our recursive methods are underpinned by the following pair of lemmas. Our proof of both lemmas uses the Dixmier long exact sequence in the setting of the Dani--Mainkar construction, as studied by Pouseele and Tirao \cite{PouseeleTirao09} and the first and second author \cite{AldiBevins24}.

\begin{lem}\label{lem:10} For all $n\ge 2$ and $k\ge 0$,
\begin{equation}
H_n^k\cong I_n(H^k_{n-1})\oplus A_n(H_{n-1}^{k-1})\oplus B_n(H_{n-1}^{k-2})\oplus C_n(\kappa_{n-1}^{k-1}). 
\end{equation}
\end{lem}

\begin{proof}
Consider the ideal of $\mathcal L(P_n)$ generated by 
$x_1,\ldots,x_{n-1},x_{1,2},\ldots,x_{n-1,n}$. The associated Dixmier long exact sequence has the form 
\begin{equation}\label{eq:8bis}
\cdots\rightarrow H^{k-1}(\mathcal D^\bullet)\xrightarrow{\delta^{k-1}}H^{k-1}(\mathcal D^\bullet)\xrightarrow{x_n^*\cdot}H^k(\mathcal L(P_n))\xrightarrow{\pi} H^k(\mathcal D^\bullet)\xrightarrow{\delta^k} H^{k}(\mathcal D^\bullet)\rightarrow \cdots,
\end{equation}
where
\begin{enumerate}[(I)]
\item $\mathcal D^\bullet\cong \mathcal C^\bullet(P_{n-1})\oplus x_{n-1,n}^*\mathcal C^\bullet(P_{n-1})$ with differential $D$ such that
\begin{equation}\label{eq:9bis}
D(\theta+x_{n-1,n}^*\omega)=Q(\theta)-x_{n-1,n}^*Q(\omega)
\end{equation}
for every $\theta,\omega\in \mathcal  C^\bullet(P_{n-1})$;
\item $x_{n}^*\cdot$ denotes left multiplication by $x_{n}^*$;
\item $\pi$ is induced by the natural projection $\frac{\partial}{\partial x_{n}^*}\circ x_{n}^*\cdot \colon \mathcal C^\bullet(P_n)\to \mathcal D^\bullet$; and
\item the connecting homomorphism $\delta^\bullet$ can be written as $\frac{\partial}{\partial x_{n}^*} \circ Q$.
\end{enumerate}
Hence, $H_n^k=H^k(\mathcal L(P_n))$ is isomorphic to the direct sum of the cokernel of $\delta^{k-1}$ and the kernel of $\delta^k$. Owing to (I), upon multiplying by $x_n^*$, ${\rm coker}(\delta^{k-1})$ is naturally identified with the subspace $A_n(H_{n-1}^{k-1})\oplus B_n(H_{n-1}^{k-2})$ of $H^k_n$. Similarly, $\ker(\delta^k)$ is isomorphic to $H_{n-1}^k\oplus x_{n-1,n}^*\kappa_{n-1}^{k-1}$, and $H_{n-1}^k$ is naturally identified through $I_n$ with a subspace of $H_n^k$, while $x_{n-1,n}^*\kappa_{n-1}^{k-1}$ is isomorphic to $\kappa_{n-1}^{k-1}$. The result now follows since, by construction, $\pi\circ C_n$ is the identity on $\kappa_{n-1}^{k-1}$.
\end{proof}

For the statement of the second lemma, we recall the notation from Remark \ref{rem:12bis}.

\begin{lem}\label{lem:10bis}
For all $n\ge 2$ and $k\ge 0$,
\begin{equation}\label{eq:12}
\mathcal H_n^k \cong A_n(\mathcal H_{n-1}^{k-1})\oplus B_n(\mathcal H_{n-1}^{k-2}\oplus I_{n-1}(\mathcal H_{n-2}^{k-2}))\oplus C_n (\mathcal K_{n-1}^{k-1}\oplus I_{n-1}( \widetilde{\mathcal K}_{n-2}^{k-1})).
\end{equation}
\end{lem}

\begin{proof}
We have a short exact sequence of cochain complexes
\[0\to (x_n^* \mathcal E_{n-1}^\bullet, Q) \to \mathcal C^\bullet(P_n,V(P_n))\oplus \mathcal C^\bullet(P_{n-1},V(P_{n-1}))\to (\mathcal E_{n-1}^\bullet,D)\to 0.\]
As vector spaces,
\[
 \mathcal E_{n-1}^\bullet = \mathcal C^\bullet(P_{n-1},V(P_{n-1}))\oplus x_{n-1,n}^* \left( \mathcal C^\bullet(P_{n-1},V(P_{n-1}))\oplus\mathcal C^\bullet(P_{n-2},V(P_{n-2}))\right),
\]
and the differential $D$ is the restriction of \eqref{eq:9bis}. Since, up to a degree shift by one unit, $ (x_n^* \mathcal E_{n-1}^\bullet, Q)$ is isomorphic to $(\mathcal E_{n-1}^\bullet,D)$, we obtain the long exact sequence
\[\cdots\rightarrow H^{k-1}(\mathcal E^\bullet)\xrightarrow{\delta^{k-1}}H^{k-1}(\mathcal E^\bullet)\xrightarrow{x_n^*\cdot} \mathcal H^k_n \oplus \mathcal H_{n-1}^k\xrightarrow{\pi} H^k(\mathcal E^\bullet)\xrightarrow{\delta^k} H^{k}(\mathcal E^\bullet)\rightarrow \cdots\]
with maps as in \eqref{eq:8bis}. Arguing as in the proof of Lemma \ref{lem:10}, we obtain
\[
{\rm coker}(\delta^{k-1})\cong A_n(\mathcal H_{n-1}^{k-1})\oplus B_n(\mathcal H_{n-1}^{k-2}\oplus I_{n-1}(\mathcal H_{n-2}^{k-2}))
\]
and
\[
\ker(\delta^k)\cong \mathcal H_{n-1}^{k}\oplus (\mathcal K_{n-1}^{k-1}\oplus I_{n-1}( \widetilde{\mathcal K}_{n-2}^{k-1})).
\]
The result then follows from $\mathcal H^k_n \oplus \mathcal H_{n-1}^k\cong {\rm coker}(\delta^{k-1})\oplus \ker(\delta^k)$.
\end{proof}

\begin{example}\label{ex:11}
For values such that $n+k\le 5$, we obtain that $\mathcal H_1^1={\rm span}(A_1([1]))$, $\mathcal H_2^2={\rm span}(B_2([1]),C_2A_1([1]))$, $\mathcal H_2^3={\rm span}(B_2A_1([1]))$, and $\mathcal H_3^2={\rm span}(C_3A_1([1]))$.   
\end{example}

The next lemma comprises properties of our cohomology operators from Section \ref{sec:2.3}.

\begin{lem}\label{lem:merge}
For every $n$,

\begin{enumerate}[(I)]
    \item the map $I_n$ is injective;
    \item the map $B_n$ is injective;
    \item $A_n\circ A_{n-1}=0$; and
    \item $C_n(\kappa_{n-1}^{k-1})\cap \kappa_n^k = 0$
\end{enumerate}

\end{lem}

\begin{proof}
Firstly, (I) is a consequence of the identification of $\mathcal C^\bullet(P_n,S)=\mathcal C^\bullet(P_{n-1},S)$ for every subset $S\subseteq V(P_{n-1})$. Then, for (III), we compute that $A_n(A_{n-1}([\omega]))=[Q(-x_{n-1,n}^*\omega)]=0$. For (II), let $[\omega]\in H^k_{n-1}$ be such that $B_n([\omega])=[x_n^*x_{n-1,n}^*\omega]=0$, that is, $x_n^*x_{n-1,n}^*\omega=Q(\psi)$ for some $\psi\in \mathcal C^\bullet(P_n)$. Since 
\[\psi = \eta+x_n^*\eta'+x_{n-1,n}^*\eta''+x_n^*x_{n-1,n}^*\eta'''\]
for some $\eta,\eta',\eta'',\eta'''\in \mathcal C^\bullet(P_{n-1})$, we obtain
\[
\omega=\frac{\partial}{\partial x_{n-1,n}^*}\frac{\partial}{\partial x_n^*} Q(\psi) = Q(\eta''')
\]
and the result follows. Lastly, for (IV), let $[\omega]$ be an element of $\kappa_{n-1}^{k-1}$. It follows from Remark \ref{rem:12} that $C_n([\omega])=[x_n^*\eta+x_{n-1,n}^*\omega]$ for some $\eta\in \mathcal C^\bullet(P_{n-1})$. If $C_n(\omega)\in \kappa_n^k$, then 
\begin{equation}\label{eq:20}
[x_n^*(x_n^*\eta+x_{n-1,n}^*\omega)]=B_n([\omega])=0.
\end{equation}
By (2), $\omega$ is a coboundary, as required.
\end{proof}

We prove another pair of auxiliary lemmas for our recursive result.

\begin{lem}\label{lem:12}
The map $A_n\circ C_{n-1}$ is monic and $A_n(C_{n-1}(\kappa_{n-2}))\cap A_n(B_{n-1}(H_{n-2}^\bullet))=0$. 
\end{lem}

\begin{proof}
If $\omega\in \kappa_{n-2}^{k-2}$, then 
\begin{equation}\label{eq:13}
A_n(C_{n-1}([\omega])) = [x_n^*x_{n-1}^*\xi+x_n^*x_{n-2,n-1}^*\omega]
\end{equation}
for some $\xi\in \mathcal C^\bullet(P_{n-2})[x_{n-2,n-1}^*]$. Suppose $A_n(C_{n-1}([\omega]))$ is a coboundary, that is, 
\begin{equation}\label{eq:14bis}
x_n^*x_{n-1}^*\xi+x_n^*x_{n-2,n-1}^*\omega=Q(\psi)
\end{equation}
 for some $\psi\in \mathcal C^\bullet(P_n)$. Since the left-hand side of \eqref{eq:14bis} is divisible by $x_n^*$ and the index $n-1$ is represented in both terms, we find that $\psi$ must be a skew-commutative polynomial in $x_n^*,x_{n-1,n}^*,x_{n-1}^*,x_{n-2,n-1}^*$ with coefficients in $\mathcal C^\bullet(P_{n-2})$, that is,
\[\psi=x_n^*x_{n-2,n-1}^*\eta+x_{n-1,n}^*\eta'+x_{n-1}^*x_n^*\eta''\]
for some $\eta\in \mathcal C^\bullet(P_{n-2})$, $\eta'\in \mathcal C^\bullet(P_{n-2})[x_n^*,x_{n-1}^*,x_{n-2,n-1}^*]$ and $\eta''\in \mathcal C^\bullet(P_{n-2})[x_{n-2,n-1}^*]$. Hence,
\begin{equation}\label{eq:15}
Q(\psi)=x_n^*x_{n-1}^*x_{n-2}^*\eta+x_n^* x_{n-2,n-1}^*Q(\eta)+x_{n-1}^*x_n^*\eta'-x_{n-1,n}^*Q(\eta')+x_{n-1}^*x_n^*Q(\eta'').
\end{equation}
Matching the coefficients of \eqref{eq:13} and \eqref{eq:15}, we obtain that $\omega=Q(\eta)$. Thus, $\omega$ is cohomologous to $0$, and $A_n\circ C_{n-1}$ is injective.

To prove the second statement, suppose $A_n(C_{n-1}([\omega]))=A_n(B_{n-1}([\omega']))$. Then,
\begin{equation}\label{eq:18bis}
x_n^*(x_{n-1}^*\eta+x_{n-2,n-1}^*\omega)-x_n^*x_{n-1}^*x_{n-2,n-1}^*\omega'=x_n^*x_{n-1}^*(\eta-x_{n-2,n-1}^*\omega')+x_n^*x_{n-2,n-1}^*\omega
\end{equation}
is a coboundary. Identifying the right-hand side of \eqref{eq:18bis} with the left-hand side of \eqref{eq:14bis} and running the previous argument yields that $\omega$ is a coboundary, as required.
\end{proof}

\begin{lem}\label{lem:14} The maps
\begin{enumerate}[(I)]
\item $A_{n+1}\circ B_n\circ A_{n-1}$;
\item $A_{n+1}\circ B_n \circ B_{n-1}$; and
\item $A_{n+1}\circ B_n\circ A_{n-2}$
\end{enumerate}
vanish identically.
\end{lem}

\begin{proof} Let $\omega\in \mathcal C^\bullet (P_{n-2})$ be a cocycle. Then,
\[A_{n+1}( B_n(A_{n-1}([\omega])))=[x_{n+1}^*x_n^*x_{n-1,n}^*x_{n-1}^*\omega]=[Q(-x_{n,n+1}^*x_{n-1,n}^*x_{n-1}^*\omega)]\]
and
\[A_{n+1}(B_n(B_{n-1}([\omega])))=[x_{n+1}^*x_n^*x_{n-1,n}^*x_{n-1}^*x_{n-2,n-1}^*\omega] =[ Q(-x_{n,n+1}^*x_{n-1,n}^*x_{n-1}^*x_{n-2,n-1}^*\omega)]\]
so we obtain (I) and (II). Then, (III) is proved in a similar fashion by noticing that, if $\omega\in \mathcal C^\bullet (P_{n-3})$ is a cocycle, then
\begin{align*}
A_{n+1}(B_n(A_{n-2}([\omega])))&=[x_{n+1}^*x_n^*x_{n-1,n}^*x_{n-2}^*\omega]\\
&=[Q((-x_{n,n+1}^*x_{n-1,n}^*x_{n-2}^*-x_{n,n+1}^*x_n^*x_{n-2,n-1}^*)\omega)].
\end{align*}
\end{proof}

We also need a lemma similar to Lemma \ref{lem:12}, which rests on the following observation. 

\begin{lem}\label{lem:22}
For every $[\omega]\in H_{n-1}^\bullet$, $A_{n+1}\circ B_n([\omega])=D_{n+1}\circ I_n([x_{n-1}^*\omega])$.
\end{lem}

\begin{proof}
The claim follows from the identity
\[
x_{n+1}^*x_n^*x_{n-1,n}^*\omega = x_{n,n+1}^*x_n^*x_{n-1}^*\omega+Q(x_{n-1,n}^* x_{n,n+1}^*\omega).
\]
\end{proof}

\begin{lem}\label{lem:16}
The map $D_n\circ I_{n-1}$ is monic, and $A_{n+1}\circ B_n\circ I_{n-1}\circ B_{n-2}=D_{n+1}\circ I_n\circ A_{n-1}\circ B_{n-2}$. 
\end{lem}

\begin{proof} For the first statement, let $[\omega] \in H_{n-2}^\bullet$, and suppose that $x_{n-1,n}^*x_{n-1}^*\omega=Q(\tau)$. Expanding in the $x_n^*$, $x_{n-1,n}^*$, $x_{n-1}^*$, and $x_{n-2,n-1}^*$ variables and imposing that the indices $n$ and $n-1$ are both represented with respective multiplicities 1 and 2, we obtain
\[
\tau=x_{n-1}^*x_{n-1,n}^* \xi + x_{n-1}^*x_{n-2,n-1}^*x_n^* \xi' + x_{n-1,n}^*x_{n-2,n-1}^* \xi''.
\]
Since $Q(\tau)$ has no $x_{n-1}^*x_{n-2,n-1}^*x_n^*$ term, it must be that $Q(-\xi')=\xi''$, and thus $\omega=Q(-\xi+ x_{n-2}^*\xi')$ is a coboundary. The second statement follows from Lemma \ref{lem:22}.
\end{proof}

The last lemma that we need for our recursive formula is as follows.

\begin{lem}\label{lem:15} The maps
\begin{enumerate}[(I)]
\item $A_{n+1}\circ B_n\circ C_{n-1}$ and
\item $A_{n+1}\circ B_n\circ C_{n-2}$ 
\end{enumerate}
are injective.
\end{lem}

\begin{proof} To prove (I), let $[\omega]\in\kappa_{n-2}$ be such that $A_{n+1}(B_n(C_{n-1}([\omega])))=0$. By Lemma \ref{lem:22} and \eqref{eq:20},
\[
A_{n+1}(B_n(C_{n-1}([\omega])))=D_{n+1}(I_n(B_{n-1}([\omega]))).
\]
The result follows from Lemmas \ref{lem:16} and \ref{lem:merge}(II). For (II), let $[\omega]\in\kappa_{n-3}$. By Lemma \ref{lem:22},
\[
A_{n+1}(B_n(C_{n-2}([\omega]))) = D_{n+1}(I_n(A_{n-1}(C_{n-2}([\omega])))),
\]
which, in conjunction with Lemmas \ref{lem:12} and \ref{lem:16}, yields the claim.
\end{proof}

\begin{rem}
The statements of Lemmas \ref{lem:merge}(IV), \ref{lem:12}, and \ref{lem:15} are independent of the choices made in the definition of $C_n$ by $A_n(H_{n-1}^{k-1})\subseteq \kappa_n^k$ and Lemmas \ref{lem:merge}(III) and \ref{lem:14}, respectively.
\end{rem}

We are ready to establish the main result of this section.

\begin{theorem}\label{recursionprop}
Define the functions $r_{n,k},p_{n,k},q_{n,k},t_{n,k} \colon \mathbb N^2\to \mathbb N$ by 
\begin{align}
p_{n,k} & = r_{n-1,k-1}+t_{n-1,k-1}\label{eq:p};\\ 
q_{n,k} & = p_{n-1,k-2}+p_{n-2,k-2}+q_{n-1,k-2}+q_{n-2,k-2}+r_{n-1,k-2}\label{eq:q};\\
r_{n,k} & = r_{n-2,k-2} + t_{n-1,k-2} + t_{n-2,k-2}\label{eq:r};\\
t_{n,k} & = p_{n-1,k-1} + p_{n-2,k-1} + q_{n-1,k-1} + q_{n-2,k-1}+ r_{n-1,k-1}\label{eq:t}
\end{align}
for all $n,k\ge 3$, together with $p_{n,1}=\delta_{n,1}$, $r_{n,2}=\delta_{n,2}$, $t_{n,2}=\delta_{n,2}+\delta_{n,3}$, and $p_{n,2}=q_{n,1}=q_{n,2}=r_{n,1}=t_{n,1}=0$ for all $n\ge 1$, as well as $p_{1,k}=q_{1,k}=r_{1,k}=t_{1,k}=r_{2,k}=p_{2,k}=t_{2,k}=0$ and $q_{2,k}=\delta_{k,3}$ for all $k>2$. Then, for all $n,k\ge 1$,
\[\beta_k(P_n)=p_{n,k}+q_{n,k}+r_{n,k}+t_{n,k}.\]
\end{theorem}

\begin{proof} The initial conditions follow from Example \ref{ex:11} and the fact that $\mathcal H_n^k=0$ when $k\le 2$ and $n> 3$. Let
\begin{align*}
p_{n,k} & =\dim (A_n(\mathcal H^{k-1}_{n-1}));\\
q_{n,k} & =\dim (B_n(\mathcal H_{n-1}^{k-2}\oplus I_{n-1}( \mathcal H_{n-2}^{k-2}))\cap \ker(A_{n+1})); \\
r_{n,k} & = \dim (A_{n+1}(B_n(\mathcal H_{n-1}^{k-2}\oplus I_{n-1}(\mathcal H_{n-2}^{k-2}))));\\
t_{n,k} & = \dim (C_n(\mathcal K_{n-1}^{k-1}\oplus I_{n-1}(\widetilde{\mathcal K}_{n-2}^{k-1}))).
\end{align*}
Combining Lemmas \ref{lem:merge}(II) and \ref{lem:12}, we obtain
\[A_n(\mathcal H_{n-1}^{k-1})\cong A_n(B_{n-1}(\mathcal H_{n-2}^{k-3}\oplus I_{n-2}(\mathcal H_{n-3}^{k-3})))\oplus A_n(C_{n-1}(\mathcal K_{n-2}^{k-2}\oplus I_{n-2}(\widetilde{\mathcal K}_{n-3}^{k-2}))),\]
which, upon using Lemma \ref{lem:12}, implies \eqref{eq:p}. It follows from \eqref{eq:12} and Lemma \ref{lem:merge}(IV) that
\begin{equation}\label{eq:23}
\mathcal K_n^k \cong A_n(\mathcal H_{n-1}^{k-1})\oplus B_n(\mathcal H_{n-1}^{k-2}\oplus I_{n-1}(\mathcal H_{n-2}^{k-2})). 
\end{equation}
Similarly, using \eqref{eq:12} and Lemma \ref{lem:12}, we obtain
\begin{equation}\label{eq:24}
\widetilde{\mathcal K}_n^k \cong A_n(\mathcal H_{n-1}^{k-1})\oplus (B_n(\mathcal H_{n-1}^{k-2}\oplus I_{n-1}(\mathcal H_{n-2}^{k-2}))\cap \ker(A_{n+1})).
\end{equation}
Combining \eqref{eq:23} and \eqref{eq:24} while keeping in mind that $C_n$ is injective yields \eqref{eq:t}. Using \eqref{eq:12} and Lemma \ref{lem:14}, we obtain
\begin{align*}
A_{n+1}(B_n(\mathcal H_{n-1}^{k-2}\oplus I_{n-1}(\mathcal H_{n-2}^{k-2}) ))& \cong A_{n+1}(B_n(B_{n-2}(\mathcal H_{n-3}^{k-4}\oplus I_{n-3}(\mathcal H_{n-4}^{k-4}))))\\
&+ A_{n+1}(B_n(C_{n-1}(\mathcal K_{n-2}^{k-3}\oplus I_{n-2}(\widetilde {\mathcal K}_{n-3}^{k-3}))))\\
&+ A_{n+1}(B_n(C_{n-2}(\mathcal K_{n-3}^{k-3}\oplus I_{n-3}(\widetilde {\mathcal K}_{n-4}^{k-3})))).
\end{align*}
Suppose that 
\[0=A_{n+1}(B_n(B_{n-2}([\omega])))+A_{n+1}(B_n(C_{n-2}([\omega'])))+A_{n+1}(B_n(C_{n-1}([\omega''])))\]
for some $\omega,\omega'\in \mathcal C^\bullet(P_{n-3})$ and $\omega''\in \mathcal C^\bullet(P_{n-2})$, that is,
\[
x_{n+1}^*x_n^*x_{n-1,n}^*(x_{n-2}^*x_{n-3,n-2}^*\omega+x_{n-2}^*\eta'+x_{n-3,n-2}^*\omega'+x_{n-1}^*\eta''+x_{n-2,n-1}^*\omega'')=Q(\psi),
\]
where $\eta'\in \mathcal C^\bullet(P_{n-3})$, $\eta''\in \mathcal C^\bullet(P_{n-2})$, and $\psi\in \mathcal C^\bullet(P_{n+1})$. Arguing as in the proof of Lemma \ref{lem:15}, we see that $x_{n+1}^*x_n^*x_{n-1,n}^*(x_{n-1}^*\eta''+x_{n-2,n-1}^*\omega'')$ is a coboundary because none of the other three terms feature the index $n-1$ with multiplicity 2. Lemma \ref{lem:15} then implies that $\omega''$ is a coboundary.  Comparing the remaining terms, the first is the only one in which the index $n-2$ has multiplicity 2 and so, once $\omega''$ is set to 0, it must be a coboundary. Applying Lemma \ref{lem:15}(II), we obtain 
\begin{align*}
A_{n+1}(B_n(\mathcal H_{n-1}^{k-2}\oplus I_{n-1}(\mathcal H_{n-2}^{k-2}) ))& \cong A_{n+1}(B_n(B_{n-2}(\mathcal H_{n-3}^{k-4}\oplus I_{n-3}(\mathcal H_{n-4}^{k-4}))))\\
&\oplus A_{n+1}(B_n(C_{n-1}(\mathcal K_{n-2}^{k-3}\oplus I_{n-2}(\widetilde {\mathcal K}_{n-3}^{k-3}))))\\
&\oplus A_{n+1}(B_n(C_{n-2}(\mathcal K_{n-3}^{k-3}\oplus I_{n-3}(\widetilde {\mathcal K}_{n-4}^{k-3})))).
\end{align*}
By Lemma \ref{lem:15}, the second and third terms contribute $t_{n-1,k-2}$ and $t_{n-2,k-2}$ to $r_{n,k}$, respectively. By Lemma \ref{lem:16}, the first term contributes $r_{n-2,k-2}$ to the right-hand side of \eqref{eq:r}, which shows that \eqref{eq:r} holds. Combining Lemma \ref{lem:merge}(III) with \eqref{eq:12}, we have that $q_{n,k}+r_{n,k}$ is equal to the sum of the right-hand sides of \eqref{eq:q} and \eqref{eq:r}, so \eqref{eq:q} holds.
\end{proof}

\begin{example}
Theorem \ref{recursionprop} can be used to tabulate the values of $p_{n,k}$, $q_{n,k}$, $r_{n,k}$, $t_{n,k}$, and thus those of $\beta_k(P_n)$.

\begin{table}[H]
    \centering
    \label{tab:pnk_values}
    \setlength{\tabcolsep}{6pt}
    \renewcommand{\arraystretch}{1.2}
    \begin{tabular}{c|rrrrrrrrrrr}
        \toprule
        \boldmath$n \setminus k$ & \textbf{1} & \textbf{2} & \textbf{3} & \textbf{4} & \textbf{5} & \textbf{6} & \textbf{7} & \textbf{8} & \textbf{9} & \textbf{10} & \textbf{11} \\
        \midrule
        \textbf{1} & 1 & 0 & 0 & 0 & 0 & 0 & 0 & 0 & 0 & 0 & 0 \\
        \textbf{2} & 0 & 0 & 0 & 0 & 0 & 0 & 0 & 0 & 0 & 0 & 0 \\
        \textbf{3} & 0 & 0 & 2 & 0 & 0 & 0 & 0 & 0 & 0 & 0 & 0 \\
        \textbf{4} & 0 & 0 & 1 & 1 & 2 & 0 & 0 & 0 & 0 & 0 & 0 \\
        \textbf{5} & 0 & 0 & 0 & 0 & 7 & 3 & 2 & 0 & 0 & 0 & 0 \\
        \textbf{6} & 0 & 0 & 0 & 0 & 5 & 6 & 14 & 5 & 2 & 0 & 0 \\
        \bottomrule
    \end{tabular}
    \caption{The values of $p_{n,k}$ for $n=1,\ldots,6$ and $k=1,\ldots,11$.}
\end{table}

\begin{table}[H]
    \centering
    \label{tab:rnk_values}
    \setlength{\tabcolsep}{6pt}
    \renewcommand{\arraystretch}{1.2}
    \begin{tabular}{c|rrrrrrrrrrr}
        \toprule
        \boldmath$n \setminus k$ & \textbf{1} & \textbf{2} & \textbf{3} & \textbf{4} & \textbf{5} & \textbf{6} & \textbf{7} & \textbf{8} & \textbf{9} & \textbf{10} & \textbf{11} \\
        \midrule
        \textbf{1} & 0 & 0 & 0 & 0 & 0 & 0 & 0 & 0 & 0 & 0 & 0 \\
        \textbf{2} & 0 & 1 & 0 & 0 & 0 & 0 & 0 & 0 & 0 & 0 & 0 \\
        \textbf{3} & 0 & 0 & 0 & 1 & 0 & 0 & 0 & 0 & 0 & 0 & 0 \\
        \textbf{4} & 0 & 0 & 0 & 3 & 1 & 1 & 0 & 0 & 0 & 0 & 0 \\
        \textbf{5} & 0 & 0 & 0 & 1 & 1 & 6 & 2 & 1 & 0 & 0 & 0 \\
        \textbf{6} & 0 & 0 & 0 & 0 & 0 & 11 & 8 & 10 & 3 & 1 & 0 \\
        \bottomrule
    \end{tabular}
    \caption{The values of $r_{n,k}$ for $n=1,\ldots,6$ and $k=1,\ldots,11$.}
\end{table}

\begin{table}[H]
    \centering
    \label{tab:qnk_values}
    \setlength{\tabcolsep}{6pt}
    \renewcommand{\arraystretch}{1.2}
    \begin{tabular}{c|rrrrrrrrrrr}
        \toprule
        \boldmath$n \setminus k$ & \textbf{1} & \textbf{2} & \textbf{3} & \textbf{4} & \textbf{5} & \textbf{6} & \textbf{7} & \textbf{8} & \textbf{9} & \textbf{10} & \textbf{11} \\
        \midrule
        \textbf{1} & 0 & 0 & 0 & 0 & 0 & 0 & 0 & 0 & 0 & 0 & 0 \\
        \textbf{2} & 0 & 0 & 1 & 0 & 0 & 0 & 0 & 0 & 0 & 0 & 0 \\
        \textbf{3} & 0 & 0 & 1 & 1 & 1 & 0 & 0 & 0 & 0 & 0 & 0 \\
        \textbf{4} & 0 & 0 & 0 & 0 & 4 & 2 & 1 & 0 & 0 & 0 & 0 \\
        \textbf{5} & 0 & 0 & 0 & 0 & 4 & 5 & 8 & 3 & 1 & 0 & 0 \\
        \textbf{6} & 0 & 0 & 0 & 0 & 1 & 2 & 18 & 16 & 13 & 4 & 1 \\
        \bottomrule
    \end{tabular}
    \caption{The values of $q_{n,k}$ for $n=1,\ldots,6$ and $k=1,\ldots,11$.}
\end{table}

\begin{table}[H]
    \centering
    \label{tab:tnk_values}
    \setlength{\tabcolsep}{6pt}
    \renewcommand{\arraystretch}{1.2}
    \begin{tabular}{c|rrrrrrrrrrr}
        \toprule
        \boldmath$n \setminus k$ & \textbf{1} & \textbf{2} & \textbf{3} & \textbf{4} & \textbf{5} & \textbf{6} & \textbf{7} & \textbf{8} & \textbf{9} & \textbf{10} & \textbf{11} \\
        \midrule
        \textbf{1} & 0 & 0 & 0 & 0 & 0 & 0 & 0 & 0 & 0 & 0 & 0 \\
        \textbf{2} & 0 & 1 & 0 & 0 & 0 & 0 & 0 & 0 & 0 & 0 & 0 \\
        \textbf{3} & 0 & 1 & 1 & 1 & 0 & 0 & 0 & 0 & 0 & 0 & 0 \\
        \textbf{4} & 0 & 0 & 0 & 4 & 2 & 1 & 0 & 0 & 0 & 0 & 0 \\
        \textbf{5} & 0 & 0 & 0 & 4 & 5 & 8 & 3 & 1 & 0 & 0 & 0 \\
        \textbf{6} & 0 & 0 & 0 & 1 & 2 & 18 & 16 & 13 & 4 & 1 & 0 \\
        \bottomrule
    \end{tabular}
    \caption{The values of $t_{n,k}$ for $n=1,\ldots,6$ and $k=1,\ldots,11$.}
\end{table}

\begin{table}[H]
    \centering
    \label{tab:bnk_values}
    \setlength{\tabcolsep}{6pt}
    \renewcommand{\arraystretch}{1.2}
    \begin{tabular}{c|rrrrrrrrrrr}
        \toprule
        \boldmath$n \setminus k$ & \textbf{1} & \textbf{2} & \textbf{3} & \textbf{4} & \textbf{5} & \textbf{6} & \textbf{7} & \textbf{8} & \textbf{9} & \textbf{10} & \textbf{11} \\
        \midrule
        \textbf{1} & 1 & 0 & 0 & 0 & 0 & 0 & 0 & 0 & 0 & 0 & 0 \\
        \textbf{2} & 0 & 2 & 1 & 0 & 0 & 0 & 0 & 0 & 0 & 0 & 0 \\
        \textbf{3} & 0 & 1 & 4 & 3 & 1 & 0 & 0 & 0 & 0 & 0 & 0 \\
        \textbf{4} & 0 & 0 & 1 & 8 & 9 & 4 & 1 & 0 & 0 & 0 & 0 \\
        \textbf{5} & 0 & 0 & 0 & 5 & 17 & 22 & 15 & 5 & 1 & 0 & 0 \\
        \textbf{6} & 0 & 0 & 0 & 1 & 8 & 37 & 56 & 44 & 22 & 6 & 1 \\
        \bottomrule
    \end{tabular}
    \caption{The values of $\beta_k(P_n)$ for $n=1,\ldots,6$ and $k=1,\ldots,11$.}
\end{table}

\end{example}

\section{An explicit formula for the essential Poincar\'e polynomial of $\mathcal N(P_n)$}\label{sec:4}
We use our recursive formula for the essential Betti numbers of path graphs to obtain a recursion for the essential Poincar\'e polynomial of $\mathcal N(P_n)$. Then, we solve this recursion to produce an explicit formula for the essential Poincar\'e polynomial of $\mathcal N(P_n)$.

Firstly, we use Theorem \ref{recursionprop} to find a recurrence for the essential Betti numbers $\beta_{n,k}$. 

\begin{lem}\label{lem:28}
    The essential Betti numbers $\beta_{n,k}$ satisfy the following recurrence relation for $n,k>3$: \[\beta_{n,k} = \beta_{n-1,k-2} + 3\beta_{n-2,k-2} + \beta_{n-2,k-3} + \beta_{n-3,k-2} + \beta_{n-3,k-3}.\]
\end{lem}

\begin{proof}
    Looking at Theorem \ref{recursionprop}, we can use the fact that the shifts $a \colon x_{n,k} \mapsto x_{n-1,k}$ and $b \colon x_{n,k} \mapsto x_{n,k-1}$ commute with each other to rewrite the recurrence relations as 
    \begin{align*}
    p_{n,k} &= ab(r_{n,k}+t_{n,k});\\
    q_{n,k} &= ab^2(p_{n,k}+a\cdot p_{n,k} + q_{n,k} + a\cdot q_{n,k}+ r_{n,k});\\
    r_{n,k} &= ab^2(a \cdot r_{n,k} + t_{n,k} + a \cdot t_{n,k});\\
    t_{n,k} &= ab (p_{n,k}+a \cdot p_{n,k} + q_{n,k} + a\cdot q_{n,k} + r_{n,k}).
    \end{align*}
    It follows that $q_{n,k} = b \cdot t_{n,k}$, and thus 
    \begin{align*}
    t_{n,k}  &= ab \cdot r_{n,k} + ab(1+a)(ab(r_{n,k}+t_{n,k}) + b\cdot t_{n,k}) \\
    &= ab \cdot r_{n,k} + (1+a)(ab^2(a \cdot r_{n,k} + a \cdot t_{n,k} + t_{n,k})) \\
    &= (ab+a + 1)(r_{n,k})\\ 
    &= T(a,b) \cdot r_{n,k},
    \end{align*} 
    where $T(a,b)=ab+a+1$. Similarly, setting $Q(a,b)=bT(a,b)$, $P(a,b)=ab(1+T(a,b))$, and $R(a,b)=ab^2(1+a(3+b)+a^2(1+b))$, we get that 
    \begin{align*}
q_{n,k} & = Q(a,b)\cdot r_{n,k};\\
p_{n,k} & = P(a,b) \cdot r_{n,k};\\
r_{n,k} & = ab^2(a\cdot r_{n,k}+T(a,b)\cdot r_{n,k}+aT(a,b)\cdot r_{n,k})\\
& = R(a,b) \cdot r_{n,k}.
    \end{align*}
    Then,
\begin{align*}    
    \beta_{n,k} &= p_{n,k}+q_{n,k}+r_{n,k}+t_{n,k}\\
    &= (P+Q+R+T)\cdot r_{n,k} \\
    &= (P+Q+R+T)\cdot (R\cdot r_{n,k}) \\
    &= R\cdot (P+Q+R+T)\cdot r_{n,k} \\
    &= R \cdot \beta_{n,k}.
\end{align*}
This algebraic reasoning is valid for $n,k > 4$ since the recurrence relation of Theorem \ref{recursionprop} is valid for $n,k >2$, so unwinding $R(a,b)$ in terms of shifts gives our result for $n,k> 3$ with the following cases treated separately:
\begin{align*}
\beta_{4,4}&=8= 1 + 6 + 0 + 0 + 1=\beta_{3,2} + 3\beta_{2,2} + \beta_{2,1} + \beta_{1,2} + \beta_{1,1}\,; \\
\beta_{4,5}&=9=4 + 3 + 2 + 0 + 0=\beta_{3,3}+3\beta_{2,3}+ \beta_{2,2}+\beta_{1,3}+\beta_{1,2}\,;\\
\beta_{5,4}&=5=0+3 + 0+2 +0=\beta_{4,2}+3 \beta_{3,2}+\beta_{3,1}+\beta_{2,2}+\beta_{2,1}\,;\\
\beta_{5,5}&= 17= 1+12+1+1+2=\beta_{4,3}+3\beta_{3,3}+\beta_{3,2}+\beta_{2,3}+\beta_{2,2}\,.
\end{align*}
\end{proof}
Lemma \ref{lem:28} yields a recursion for the essential Poincar{\'e} polynomials of path graphs.

\begin{lem}\label{PoinRec}
    The essential Poincar\'e polynomial $\beta_n(x)=\beta(P_n,x)$ satisfies the following recursion relation for $n>3$:
        \begin{equation}\label{eq:38}
            \beta_{n}(x) = x^2\beta_{n-1}(x) + x^2(3+x) \beta_{n-2}(x) + x^2(1+x)\beta_{n-3}(x).
        \end{equation}
\end{lem}
\begin{proof}
    We have 
    \[
    \beta_{n,k} = ab^2(1+a(3+b) + a^2(1+b))\beta_{n,k}.
    \]
    Noting that $b \leftrightarrow x\cdot (-)$ and $a\cdot \beta_n(x)=\beta_{n-1}(x)$ and summing over $k$, we obtain 
    \begin{align*}
    \beta_n(x) & =  x^2(a+a^2(3+x) +a^3(1+x)) \beta_n(x)\\
    & = x^2\beta_{n-1}(x) + x^2(3+x) \beta_{n-2}(x) + x^2(1+x)\beta_{n-3}(x).
    \end{align*}
\end{proof}

We are now ready to solve our recurrence from Lemma \ref{PoinRec} in an explicit fashion.

\begin{theorem}\label{PoinFormula}
    For $\ell \geq 0$, define \begin{equation}\label{ThetaPoly} \Theta_\ell(x) =\sum_{\substack{0\leq k\leq m\leq n \\ n+m+k = \ell}} \binom{n}{m} \binom{m}{k}x^{2n}(3+x)^{m-k} (1+x)^k. \end{equation} If $N >2$ is an integer, then 
    \begin{equation}\label{eq:essentialexplicit}
        \beta_{N+1}(x) = x\Theta_N(x) + 2x^2 \Theta_{N-1}(x) + x^2(1+x)\Theta_{N-2}(x).   
    \end{equation}
    Furthermore,
    \begin{equation}\label{eq:43}
\sum_{n=0}^\infty \beta_{n+1}(x)t^n = \frac{x+2x^2t+x^2(1+x)t^2}{1-x^2t-x^2(3+x)t^2-x^2(1+x)t^3}.
    \end{equation}
\end{theorem}

\begin{proof}
    Lemma \ref{PoinRec} yields that \begin{equation}\label{brec} \beta_{n+4}(x) = x^2 \beta_{n+3}(x) + x^2(3+x)\beta_{n+2}(x)+x^2(1+x)\beta_{n+1}(x).
    \end{equation}  
    Consider now the formal power series in the formal variable $t$ with coefficients in $\mathbb{C}[x]$
    \[
    g(x,t)= \sum_{n=0}^\infty \beta_{n+1}(x) t^n .
    \]
    Multiplying both sides of \eqref{brec} by $t^{n+3}$, we obtain 
    \[   
    \beta_{n+4}(x) t^{n+3} = x^2 t \cdot \beta_{n+3}(x)t^{n+2} + x^2(3+x) t^2\cdot \beta_{n+2}(x)t^{n+1} + x^2(1+x)t^3 \cdot \beta_{n+1}(x) t^n.
    \]
    Summing from $n=0$ to $\infty$, we find that 
    \[ g(x,t) - \beta_1(x) - \beta_2(x) t - \beta_3(x) t^2
    \]
    is equal to
    \[
    x^2t (g(x,t)-\beta_1(x) - \beta_2(x) t) + x^2(3+x) t^2 (g(x,t)-\beta_1(x)) + x^2(1+x)t^3g(x,t).
    \] 
    Moving $g$ to the left, we see that \[(1-x^2t - x^2(3+x)t^2 - x^2(1+x)t^3)g(x,t)\]
    is equal to
    \[
\beta_1(x) + (\beta_2(x) -x^2 \beta_1(x)) t + (\beta_3(x) - x^2 \beta_2(x) - x^2(3+x)\beta_1(x))t^2.
    \]
    Recalling that $\beta_1 = x$, $\beta_2 = 2x^2+x^3$, and $\beta_3 = x^2 + 4x^3 + 3x^4+x^5$, we obtain 
    \[ 
    (1-x^2 t- x^2(3+x) t^2 - x^2(1+x)t^3)g = x + 2x^2 t + x^2(1+x)t^2.
    \]
    If we let 
    \[\nu(x,t) = x^2 t + x^2 (3+x)t^2 + x^2(1+x)t^3,
    \]
    then we can formally invert $1-\nu(x,t)$ and obtain 
    \begin{equation}\label{gEq} g(x,t)  =  \dfrac{x+2x^2 t + x^2(1+x)t^2}{1- \nu(x,t)},\end{equation} 
    from which \eqref{eq:43} follows. Now, 
    \begin{align*}
    \dfrac{1}{1-\nu(x,t)} &= \sum_{n=0 }^\infty (\nu(x,t))^n \\
    &= \sum_{n=0}^\infty x^{2n} (1+ t((3+x) + t(1+x)) )^n  t^n\\
    &= \sum_{n=0}^\infty \sum_{m=0}^n \binom{n}{m} x^{2n}  ((3+x) + t (1+x))^m t^{n+m}\\
    &= \sum_{n=0}^\infty \sum_{m=0}^n \sum_{k=0}^m \binom{n}{m} \binom{m}{k} x^{2n} (3+x)^{m-k} (1+x)^k t^{n+m+k}.
    \end{align*}
    We recover our definition of $\Theta_\ell(x)$ by looking at the coefficient of $t^\ell$. Thus, 
    \[
    \dfrac{1}{1-\nu(x,t)} = \sum_{\ell=0}^\infty \Theta_\ell(x) t^\ell.
    \] 
    Plugging into \eqref{gEq}, we get 
    \[ 
    \sum_{\ell=0}^\infty  \beta_{\ell+1}(x) t^\ell =  \sum_{\ell=0}^\infty x \Theta_\ell(x) t^\ell + \sum_{\ell=0}^\infty 2x^2 \Theta_\ell(x) t^{\ell+1} + \sum_{\ell=0}^\infty x^2(1+x)\Theta_\ell(x) t^{\ell+2}.
    \] 
    Grouping coefficients of $t^N$ finishes the proof.
\end{proof}

Using these explicit formulas, we can describe the degree and value of the lowest nonzero essential Betti numbers.
\begin{cor}\label{ordvancor}
   For every $r$, we have
   \begin{align*}
      \beta_{3r+1}(x) &= x^{2r+1} + \text{higher-order terms}; \\
      \beta_{3r+2}(x) &= (3r+2)x^{2r+2} + \text{higher-order terms}; \\
      \beta_{3r+3}(x) &= x^{2r+2}+ \text{higher-order terms}.
   \end{align*}
\end{cor}
\begin{proof}
    Observe that, in \eqref{ThetaPoly}, the triples $(n,m,k)$ contributing lowest-order terms to $\Theta_\ell(x)$ are precisely those with minimal $n$. In the case $\ell=3r$, the only triple minimizing $n$ is $(n,m,k)=(r,r,r)$, with a corresponding lowest-order term
    $\binom{r}{r}\binom{r}{r}x^{2r}3^{r-r}1^r=x^{2r}.$
    Similarly, in the case $\ell=3r+1$, we have triples $(n,m,k)=(r+1,r+1,r-1)$ and $(r+1,r,r)$, which contribute lowest-order terms $\frac{9}{2}r(r+1)x^{2r+2}$ and $(r+1)x^{2r+2}$, respectively. In the case $\ell=3r+2$, only the triple $(n,m,k)=(r+1,r+1,r)$ minimizes $n$, and it contributes the lowest-order term $3(r+1)x^{2r+2}$. Thus, we have
    \begin{align*}
        \Theta_{3r}(x) &= x^{2r} + \text{higher-order terms}; \\
        \Theta_{3r+1}(x) &= \left(\frac92r+1\right)(r+1)x^{2r+2} + \text{higher-order terms}; \\
        \Theta_{3r+2}(x) &= 3(r+1)x^{2r+2} + \text{higher-order terms}.
    \end{align*}
    Plugging these formulas into \eqref{eq:essentialexplicit} yields the claim.
\end{proof}

\begin{rem}\label{ordvanprop} We can also prove Corollary \ref{ordvancor} independently of our explicit formula for the essential Poincar\'e polynomials from Theorem \ref{PoinFormula} by using the order of vanishing of a polynomial at $0$, which we denote by $\operatorname{ord}$, instead. Applying $\operatorname{ord}$ to both sides of \eqref{eq:38}, defining $\alpha_n = \operatorname{ord}(\beta_n)$, and using the fact that $\operatorname{ord}$ is a valuation, we obtain \[\alpha_n \geq 2 + \min(\alpha_{n-1},\alpha_{n-2},\alpha_{n-3}).\] Now, we can use the fact that our polynomials have coefficients which are non-negative and the recurrence relation only involves positive integers; the cancellation of terms, which is the only reason the previous inequality can be strict, cannot occur, and thus \[\alpha_n = 2 + \min(\alpha_{n-1},\alpha_{n-2},\alpha_{n-3}).\]
    Since the sequence $\alpha_n$ begins with $\alpha_1=1, \alpha_2 = 2,\alpha_3=2$, we can use induction to prove that $|\alpha_{n+1}-\alpha_n| \leq 1$ and that $\alpha$ is a monotonically increasing sequence. Suppose this statement is true for all $k < m$ for a given $m>3$. Then, 
    \[
    \alpha_m = 2 + \min(\alpha_{m-1},\alpha_{m-2},\alpha_{m-3}) = 2 + \alpha_{m-3}.
    \] Now, \[\alpha_{m-1} = 2+ \alpha_{m-4} \geq 1 + \alpha_{m-3} = \alpha_m -1\] implies that $\alpha_m \leq \alpha_{m-1}+1$, and $\alpha_{m-1}\leq 2 + \alpha_{m-3} = \alpha_m $ Thus, $\alpha_n$ is monotonically increasing, which means that our recurrence relation is \[\alpha_{n} = 2 + \alpha_{n-3}\] with initial conditions $\alpha_1 = 1, \alpha_2=2,\alpha_3 = 2$. Letting $n = 3 \ell + \varepsilon$ for fixed $\varepsilon \in \{1,2,3\}$ and solving the equation $\alpha'_\ell = 2 + \alpha'_{\ell-1}$ for $\alpha'_\ell = \alpha_{3 \ell + \varepsilon}$ proves the first set of claims about $\operatorname{ord}(\beta_n)$. If we are looking for the minimal terms of 
    \[
    \beta_{3\ell + \varepsilon }(x)= x^2(\beta_{3(\ell-1)+2 + \varepsilon}(x) + (3+x) \beta_{3(\ell-1) + 1 + \varepsilon}(x) + (1+x)\beta_{3(\ell-1)+ \varepsilon}(x)),
    \]  
    then, since our polynomials have positive coefficients the only contributions (ignoring the factor of $x^2$ for now) to the minimal term can come from the sum 
    \[
    \beta_{3(\ell-1)+2 + \varepsilon}(x) + 3 \beta_{3( \ell-1)+1 + \varepsilon}(x) + \beta_{3(\ell-1)+\varepsilon}(x).
    \]
    If we define $ c^{1}_\ell, c^{2}_\ell, c^{3}_\ell$ as \begin{align*} 
    \beta_{3\ell+1}(x) &= c^{1}_\ell x^{2\ell+1} + \text{higher-order terms};\\
    \beta_{3\ell+2}(x) &= c^{2}_\ell x^{2\ell+2} + \text{higher-order terms};\\
    \beta_{3\ell+3}(x) & = c^{3}_\ell x^{2\ell+2} + \text{higher-order terms},
    \end{align*}
    then we get the system of equations 
    \begin{align*}
    c^{1}_\ell x^{2\ell+1} + \text{higher-order terms}&= c^3_{\ell-1} x^{2\ell+2} + 3c^{2}_{\ell-1}x^{2\ell+2} + c^{1}_{\ell-1}x^{2\ell+1}+ \text{higher-order terms};\\
    c^2_\ell x^{2\ell+2} + \text{higher-order terms} &= c^1_{\ell}x^{2\ell+3} + 3 c^3_{\ell-1} x^{2\ell+2} + c^2_{\ell-1}x^{2\ell+2} + \text{higher-order terms};\\
    c^3_\ell x^{2\ell+2} + \text{higher-order terms} &= c^2_{\ell} x^{2\ell+4} + 3c^1_\ell x^{2\ell+3} + c^3_{\ell-1} x^{2\ell+2} + \text{higher-order terms}.
    \end{align*}
    Now, equating the minimal terms, we obtain that $c^1_\ell = c^1_{\ell-1}$, $c^2_\ell = 3 c^3_{\ell-1} + c^2_{\ell-1}$, $c^3_\ell = c^3_{\ell-1}$. Using $c^1_0 = 1$, $c^2_0 = 2$, and $c^3_0 = 1$, we find that $c^1_\ell = 1$, $c^2_\ell = 3\ell+2$, and $c^3_\ell = 1$.
\end{rem}

\section{An explicit formula for the Poincar\'e polynomial of $\mathcal N(P_n)$}\label{sec:5}
In this section, we emulate our line of computation from Section \ref{sec:4} to give a recursion for the Poincar\'e polynomials of path graphs, which we proceed to solve explicitly. Our results hinge on our work for essential Poincar\'e polynomials. Moreover, we use our explicit formula to show that the Poincar\'e polynomials of path graphs are unimodal.

Firstly, we note that, by \eqref{eq:3}, the calculation of the Betti numbers of $\mathcal L(G)$ can be, at least in principle, reduced to the calculation of the essential Betti numbers. In general, this task is challenging, for it requires detailed knowledge of the structure of the induced subgraphs. However, in the case of paths, it can be carried out explicitly. 

\begin{lem}\label{rem:33}
    Define $\beta_{0}(x) = b_{0}(x)= 1$ and $b_n(x)=b(P_n,x)$ for all $n\ge 1$. Then, for all $n >1$, \begin{equation}\label{PoincareBadRecurs} b_{n}(x) = \beta_n(x)+\sum_{k=1}^{n}\beta_{n-k} (x) \cdot b_{k-1}(x).\end{equation}
\end{lem}
\begin{proof}
     If $H \leq P_{n}$ is an induced subgraph, then it either contains $x_n$ or it does not. If it does not, then $H$ is an induced subgraph of $P_{n-1}$. Thus, 
     \[
     b_{n}(x) = \sum_{H \leq P_{n-1}} \beta(H,x)+b'_n(x),
     \]
     where $b'_n(x)$ only counts the contributions of induced subgraphs containing $x_n$. We can recognize $b_{n-1}(x)$ to see that $b_n(x) = b_{n-1}(x)+b'_n(x)$.
     
     If $H$ contains $x_n$, then it contains a maximal connected component $\{x_{k+1},x_{k+2},\dots,x_{n}\}$ for some $0\leq k \leq n-1$. Fixing this connected component, we see that the possible $H$ with this connected component are the union of the connected component and an induced subgraph of $P_{k-1}$ because the vertex $x_k$ must be excluded; otherwise, $H$ would not be induced. Thus, the contribution of such $H$ to $b'_n(x)$ is 
     \[\sum_{H' \leq P_{k-1}} \beta(H',x) \cdot \beta_{n-k}(x) = \beta_{n-k}(x)\cdot b_{k-1}(x).
     \]  
     Since this collection of subsets of induced subgraphs is disjoint and exhausts all $H$ which contain $x_n$, we obtain
     \[
     b_n(x) = \beta_{n}(x) + \sum_{k=1}^n \beta_{n-k}(x) \cdot b_{k-1}(x).
     \]
\end{proof}
In addition to the recursive relation of Lemma \ref{rem:33}, we have a formula expressing the Betti numbers of $\mathcal N(P_n)$ in terms of the essential Betti numbers of $\mathcal N(P_n)$.

\begin{prop}
If $P_n$ is the path graph on $n$ vertices, then
\[b_k(P_n)= \sum_{\mathbf m} c_{\mathbf m} \sum_{\mathbf k} \beta_{k_1}(P_1)\cdots\beta_{k_{m_1}}(P_1)\beta_{k_{m_1+1}}(P_2)\cdots \beta_{k_{m_1+\cdots+m_N}}(P_N),\]
where the first sum is over all non-negative integer vectors $\mathbf m=(m_1,\ldots,m_N)$ such that $m_1+2m_2+\cdots+Nm_N\le 2k-1$, the second sum is over all positive integer vectors $\mathbf k=(k_1,\ldots,k_{m_1+\cdots+m_N})$ such that $k=k_1+\cdots+k_{m_1+\cdots+m_N}$, and
\[c_{\mathbf m}= \frac{(n-(m_1+2m_2+\cdots+Nm_N)+1)!}{(n-(2m_1+3m_2+\cdots+(N+1)m_N)+1)!m_1!m_2!\cdots m_N!}.\]
\end{prop}

\begin{proof} 
Every induced subgraph of $P_n$ is a disjoint union of path graphs, so \eqref{eq:3} reduces to
\begin{equation}\label{eq:7}
b_k(P_n)=\sum_{m_1+2m_2+\cdots+Nm_N\le 2k-1} \left|\binom{P_n}{m_1P_1+\cdots+m_NP_N}\right| \beta_k(m_1P_1+\cdots+m_NP_N).
\end{equation}
Repeated application of \eqref{eq:4} yields
\begin{equation}\label{eq:8}
\beta_k(m_1P_1+\cdots+m_NP_N)=\sum_{k_1+\cdots+k_{m_1+\cdots+m_N}=k} \beta_{k_1}(P_1)\cdots\beta_{k_{m_1}}(P_1)\beta_{k_{m_1+1}}(P_2)\cdots \beta_{k_{m_1+\cdots+m_N}}(P_N).
\end{equation}
On the other hand, by contracting each connected component of $m_1P_1+\cdots+m_NP_N$ to a single vertex, we see that 
\begin{equation}\label{eq:9}
\left|\binom{P_n}{m_1P_1+\cdots+m_NP_N}\right| = s_{m_1+\cdots+m_N}(P_{n-(m_2+2m_3+\cdots+(N-1)m_N)})\binom{m_1+\cdots+m_N}{m_1,\ldots,m_N},
\end{equation}
where $s_k(G)$ denotes the number of independent sets of cardinality $k$ in a graph $G$. The result follows by \eqref{eq:7}, \eqref{eq:8}, \eqref{eq:9}, and Arocha's \cite{Arocha84} formula $s_k(P_n)=\binom{n-k+1}{k}$.
\end{proof}

The recurrence relation of Lemma \ref{rem:33} is too unwieldy for our purposes. Nevertheless, our coefficients of each $b_{k-1}(x)$ each have their own recurrence relation that we have analyzed, which we can use to our advantage. 

\begin{prop}\label{totalPoinRecurs}
    The Poincar\'e polynomials $b_n(x)$ satisfy the following recurrence relation when $n \geq 1$:
    \begin{equation}\label{eq:53}
        b_{n+2}(x) = (1+x^2)b_{n+1}(x) + x(1+x)^2 b_n(x).
    \end{equation}
\end{prop}

\begin{proof}
Our recurrence relation for $\beta_n(x)$ reads as \[
\beta_{n+3}(x) = x^2\beta_{n+2}(x) + x^2(3+x) \beta_{n+1}(x) + x^2(1+x)\beta_n(x),
\] 
so we take \eqref{PoincareBadRecurs} and shift up by $3$ to obtain
\begin{align*}
b_{n+3}(x) & = \beta_{n+3}(x)+\sum_{k=1}^{n+3}\beta_{n-k+3}(x) b_{k-1}(x) \\
&=\beta_{n+3}(x)+ b_{n+2}(x) + \beta_1(x) b_{n+1}(x) + \beta_2(x) b_n(x) + \beta_3(x) b_{n-1}(x) + \\
&\quad +\sum_{k=1}^{n-1} \beta_{n-k+3}(x) b_{k-1}(x)\\
& =\beta_{n+3}(x)+ b_{n+2}(x) + \beta_1(x) b_{n+1}(x) + \beta_{2}(x) b_{n}(x) + \beta_{3}(x)b_{n-1}(x) + \\
&+x^2 \left(\sum_{k=1}^{n-1}(\beta_{n-k+2}(x) + (3+x)\beta_{n-k+1}(x) + (1+x)\beta_{n-k}(x))b_{k-1}(x)\right).
\end{align*}

We can separate the terms on the right-hand side into the terms 
\begin{align*}
x^2\sum_{k=1}^{n-1}\beta_{n-k+2}(x) \cdot b_{k-1}(x) &= x^2\left ( b_{n+2}(x) - \beta_{n+2}(x)-b_{n+1}(x) - \beta_{1}(x)b_{n}(x) - \beta_2(x) b_{n-1}(x) \right );\\
x^2(3+x)\sum_{k=1}^{n-1} \beta_{n-k+1}(x)\cdot b_{k-1}(x) &= x^2(3+x)(b_{n+1}(x) -\beta_{n+1}(x)- b_n(x) - \beta_1(x) b_{n-1}(x));\\
x^2(1+x) \sum_{k=1}^{n-1} \beta_{n-k}(x) \cdot b_{k-1}(x) &= x^2(1+x) (b_n(x) - \beta_n(x)-b_{n-1}(x)).
\end{align*}
Putting everything together, we obtain 
\[
b_{n+3}(x) = (1+x^2)b_{n+2}(x) + \gamma_0(x) \cdot b_{n+1}(x) + \gamma_1(x)\cdot b_n(x) + \gamma_2(x) \cdot b_{n-1}(x),
\] 
where 
\begin{align*}
\gamma_0 (x) & = \beta_1 (x)-x^2 + x^2(3+x)\\
&= \beta_1(x) + 2x^2 + x^3;\\
\gamma_1(x) &= \beta_2(x) - x^2( \beta_1(x) + (3+x) - (1+x)) \\
&= \beta_2(x)- x^2(\beta_1(x) + 2); \\
\gamma_2(x) &= \beta_3(x) - x^2(\beta_2(x) + (3+x)\beta_1(x) + (1+x)).
\end{align*}
Plugging in $\beta_1(x) = x$, $\beta_2 (x)= 2x^2 + x^3$, and $\beta_3(x) = x^2 + 4x^3+3x^4 + x^5$ gives $\gamma_1(x) = \gamma_2(x) =0 $ and $\gamma_0 (x)= x + 2x^2 + x^3 = x(1+x)^2$. Lastly, shift $n$ back by $1$.
\end{proof}

\begin{rem}
    The proof implies that $\beta_3(x)$ follows the recurrence relation proven for larger $n$ given $\beta_0(x) =1$, so we did not actually  need to sequester $\beta_3(x) b_{n-1}(x)$ at the start.
\end{rem}

\begin{table}[H]
    \centering
    \label{tab:bn_polys}
    \renewcommand{\arraystretch}{1.3}
    \begin{tabular}{c|l}
        \toprule
        \boldmath$n$ & \textbf{$b_n(x)$} \\
        \midrule
        \textbf{1} & $1+x$ \\
        \textbf{2} & $1+2x+2x^2+x^3$ \\
        \textbf{3} & $1+3x+6x^2+6x^3+3x^4+x^5$ \\
        \textbf{4} & $1+4x+11x^2+16x^3+16x^4+11x^5+4x^6+x^7$ \\
        \textbf{5} & $1+5x+17x^2+33x^3+48x^4+48x^5+33x^6+17x^7+5x^8+x^9$ \\
        \textbf{6} & $1+6x+24x^2+58x^3+107x^4+140x^5+140x^6+107x^7+58x^8+24x^9+6x^{10}+x^{11}$ \\
        \bottomrule
    \end{tabular}
    \caption{The Poincar\'e polynomials $b_n(x)$ of $\mathcal N(P_n)$ for $n=1,\ldots,6$.}
\end{table}

We take a big step toward an explicit formula for the Poincar{\'e} polynomial of $\mathcal{N}(P_{n})$.

\begin{prop}\label{PoincarePolyFormula}
    Let \[\Psi_\ell(x)=  \sum_{n=\lceil\frac{\ell}{2}\rceil }^\ell \binom{n}{2n-\ell} x^{\ell-n} (1+x^2)^{2n-\ell} (1+x)^{2(\ell-n)}.\] Then, for $n \geq 1$, 
    \[
    b_{n+1}(x) =(1+x) \left(\Psi_n(x) + x \Psi_{n-1}(x)\right). 
    \]
    Moreover,
    \begin{equation}\label{eq:53tris}
    \sum_{n=0}^\infty b_{n+1}(x)t^n = \frac{(1+x)(1+xt)}{1-(1+x^2)t-x(1+x)^2t^2}.
    \end{equation}
\end{prop}

\begin{proof}
    We use generating functions to construct a solution to the recurrence relation 
    \[
    b_{n+2}(x)= (1+x^2) b_{n+1} (x)+ x(1+x)^2b_n(x).
    \] Let 
    \[
    g(x,t)= \sum_{n=0}^\infty b_{n+1}(x) t^n.
    \]  
    Multiplying by $t^{n+2}$, we obtain 
    \[
    b_{n+3}(x)t^{n+2} = t(1+x^2)b_{n+2}(x)t^{n+1} + t^2x(1+x)^2 b_{n+1}(x) t^n.
    \] 
    Taking the formal series from $n=0$, we get 
    \[
    (g(x,t)-b_2(x) t - b_1(x)) = t(1+x^2)(g(x,t) - b_1(x)) + t^2x(1+x)^2 g(x,t).
    \] 
    We put all terms involving $g(x,t)$ on the left, and we recall that $b_1(x) = 1+x$. Then,
    \[
    b_2(x) = 1+2x + 2x^2 + x^3 = (1+x)(1+x+x^2)
    \]
    leads to 
    \[\left(1- (1+x^2)t - x(1+x)^2t^2 \right)g(x,t) = b_1(x) + (b_2(x) - (1+x^2)b_1(x))t = b_1(x)(1 + xt).
    \] 
    Let 
    \begin{equation}\label{eq:53bis}
    \mu(x,t) = (1+x^2)t + x(1+x)^2t^2 = t(1+x^2  + t x(1+x)^2).
    \end{equation}
    Then,
    \begin{align*}
    \dfrac{1}{1-\mu(x,t)} &= \sum_{m=0}^\infty (\mu(x,t))^m \\
    &= \sum_{m=0}^\infty ((1+x^2)+tx(1+x)^2)^m t^m\\
    &= \sum_{m=0}^\infty \sum_{k=0}^m \binom{m}{k} (1+x^2)^{m-k} x^k (1+x)^{2k} t^{m+k}.
    \end{align*}
    The $t^{\ell}$ coefficient of $\dfrac{1}{1-\mu(x,t)}$ is 
    
    \begin{align*}
    \Psi_\ell(x) &= [t^{\ell}]\left(\dfrac{1}{1-\mu(x,t)}\right) \\
    & = \sum_{\substack{0\leq k \leq m \\ k+m=\ell } } \binom{m}{k} x^k (1+x^2)^{m-k}(1+x)^{2k}\\
    &= \sum_{n=0}^{\ell} \binom{n}{\ell-n} x^{\ell-n}(1+x^2)^{2n-\ell} (1+x)^{2\ell - 2n} \\ 
    &= \sum_{n=\lceil\frac{\ell}{2}\rceil }^\ell \binom{n}{2n-\ell} x^{\ell-n} (1+x^2)^{2n-\ell} (1+x)^{2(\ell-n)}.
    \end{align*} 
    Putting it back into our equation for $g$, we find that
    \begin{align}
    g(x,t)  &= b_1(x)\left( \dfrac{1}{1-\mu(x,t)}  + x  \dfrac{t}{1-\mu(x,t)}\right)\label{eq:53quater} \\ 
    & =  (1+x) \left( \sum_{\ell=0}^{\infty} \Psi_\ell(x) t^\ell + \sum_{\ell=0}^\infty x \Psi_\ell(x) t^{\ell +1}\right).
    \end{align}
    We look at the coefficients of $t^n$ to find that
    \[
    b_{n+1}(x) = (1+x) \left(\Psi_n(x) + x \Psi_{n-1}(x)\right). 
    \]
  Finally, \eqref{eq:53tris} follows from combining \eqref{eq:53bis} and \eqref{eq:53quater}.
\end{proof}

We arrive at the main result of this section.

\begin{theorem}\label{CharacteristicPolynomialForm}
    Let \[\delta(x)= \sqrt{(1+x^2)^2 + 4x(1+x)^2} = \sqrt{(1+x)^4 + 4x^2}\]
    and
    \begin{equation}\label{eq:54}
    \psi_{\pm}(x) = \dfrac{1}{2}\left(1+x^2 \pm \delta(x)\right).
    \end{equation} 
   Then, for all $n\ge 1$, \[b_n(x) =\frac{1+x}{\delta(x)}\left(\big((\psi_+(x))^n-(\psi_-(x))^n\big)+x\big((\psi_+(x))^{n-1}-(\psi_-(x))^{n-1}\big)\right).\]
\end{theorem}

\begin{proof}
    By \eqref{eq:53}, our characteristic polynomial equation is 
    \begin{equation}\label{eq:55bis}
    (\psi_{\pm}(x))^2 - (1+x^2)\psi_{\pm}(x) - x(1+x)^2 = 0.
    \end{equation}
    The solution to this equation in algebraic functions is 
    \[ 
    \psi_{\pm}(x) = \dfrac{1}{2}\left(1+x^2 \pm \sqrt{(1+x^2)^2 + 4 x(1+x)^2} \right) = \dfrac{1}{2}(1+x^2 \pm \delta(x)).
    \] 
    We arrive at the general formula 
    \begin{equation}\label{eq:55}
    b_{n}(x) = A(x)\cdot (\psi_+(x))^n + B(x) \cdot (\psi_-(x))^n
    \end{equation}
    for algebraic $A(x),B(x)$. Using $b_1(x) = 1+x$ and $b_2(x) = 1+2x+2x^2+x^3$, we obtain 
    \begin{align}
    A(x)\psi_+(x) + B(x) \psi_-(x) & = 1+x;\\
    A(x) (\psi_+(x))^2 + B(x) (\psi_-(x))^2 &= 1+2x+2x^2 + x^3.\label{eq:56}
    \end{align}
    Substituting 
    \[
    B(x) (\psi_-(x))^2 = \psi_-(x)\cdot (1+x - A(x) \psi_+(x))
    \]
    into \eqref{eq:56}, we get 
    \[
    A(x)(\psi_+(x))^2 + \psi_-(x)\cdot (1+x) - A(x) \psi_+(x) \psi_-(x) = (1+x)(1+x+x^2),
    \]
    which implies
    \begin{align*}
    A(x) &= \dfrac{ (1+x)(2+2x+2x^2 - 1 - x^2 + \delta(x)) }{2\psi_+(x) (\psi_+(x) - \psi_-(x))}\\ 
    &=\dfrac{(1+x)(x + \psi_+(x)) }{\delta(x) \psi_+(x)}.
    \end{align*}
    Therefore,
    \[
    B(x) = \dfrac{1+x - A(x) \psi_+(x)}{\psi_-(x)} = \dfrac{(x+1) ( \delta(x) - x - \psi_+(x))}{\delta(x) \psi_-(x)}= -\dfrac{(x+1)(x + \psi_-(x))}{\delta(x) \psi_-(x)}.
    \] 
\end{proof}

We conclude this section with a unimodality application of Theorem \ref{CharacteristicPolynomialForm}.

\begin{prop}\label{prop:38}
    For $n \geq 1$, the Poincar\'e polynomial $b_n(x)$ is unimodal; for all $n\ge 2$ and $1\le  k < n$, $b_{n,k}>b_{n,k-1}$. 
\end{prop}

\begin{proof}Let $\Delta_{n,k}= b_{n,k} - b_{n,k-1}$ for $1\leq k<n$. We tabulate those values for $n \leq 5$.
    \begin{table}[H]
        \centering
        \begin{tabular}{|c|c|c|c|c|c|}
        \hline
             $\Delta_{n,k}$& $n=1$ & 2 & 3 & 4 & 5 \\
             \hline
             $k=1$         &       & 1 & 2 & 3 & 4 \\
             \hline
             $2$           &       &   & 3 & 7 & 12\\
             \hline
             $3$           &       &   &   & 5 & 16\\
             \hline
             $4$           &       &   &   &   & 15\\
             \hline
        \end{tabular}  
        \caption{The differences $\Delta_{n,k}$ for $n \leq 5$.}
    \end{table}
    We can construct a recursion relation for $\Delta_{n,k}$ using the one for $b_n$: \[\Delta_{n,k} = \Delta_{n-1,k} + \Delta_{n-1,k-2} + \Delta_{n-2,k-1}+2\Delta_{n-2,k-2} + \Delta_{n-2,k-3}. \] For $k\leq n-2$, we face no issue, but for $k=n-1$, we have the symbols $\Delta_{n-1,n-1}$ and $\Delta_{n-2,k-1}$ which are outside the scope of our definition for $\Delta$. By Poincar\'e duality, $\Delta_{n-1,n-1} = b_{n-1,n-1} - b_{n-1,n-2} = 0$ and $\Delta_{n-2,n-2} = b_{n-2,n-2} - b_{n-2,n-3} = 0$. We solve explicitly for the first few terms of $b_n$. 
    \begin{enumerate}[(I)]
        \item The term $b_{n,1}$ is equal to $n$: \[ b_{n,1} = b_{n-1,1} + b_{n-2,0} = b_{n-1,1}+1.\]
        \item The term $b_{n,2}$ is $ \frac{1}{2}(n^2 +3n) - 3$: \[b_{n,2} = b_{n-1,2} + 1 +b_{n-2,1} + 2 = b_{n-1,2}+ n-2 + 3 = b_{n-1,2} + n+1.\]
        
    \end{enumerate}
With these results in hand, suppose that $N$ is a positive integer larger than $4$ such that the statement $\Delta_{n,k} >0$ is true for $ n < N$ and $k \leq n-1$. Then, if $4 \leq k \leq N-2$, \[ \Delta_{N,k} = \Delta_{N-1,k} + \Delta_{N-1,k-2} + \Delta_{N-2,k-1} + 2\Delta_{N-2,k-2} + \Delta_{N-2,k-3} >0. \] For $k = N-1$, we have that \[\Delta_{N,N-1} = \Delta_{N-1,N-3} + 2 \Delta_{N-2,N-3} + \Delta_{N-2,N-4} > 0.\] For the other cases, we have 
\begin{align*}
\Delta_{N,1} &= N-1;\\
\Delta_{N,2} &= \frac{1}{2}(N^2 +N) -3 \geq 12 > 0;\\
\Delta_{N,3} &= b_{N,3} - b_{N,2} \\
& = b_{N-1,3} + b_{N-1,1} + b_{N-2,2} + 2b_{N-2,1} + 1 -(b_{N-1,2} + 1 + b_{N-2,1}+ 2)\\ 
&=\Delta_{N-1,3} +\Delta_{N-2,2} + b_{N-1,1} + 2b_{N-2,1} - 2 \\
&= \Delta_{N-1,3} + \Delta_{N-2,2} + 3N - 7 > 0. 
\end{align*}
\end{proof}

\section{The asymptotic behavior of the cohomology of path graphs}\label{sec:6}
We end this paper with some asymptotics on the total Betti number of $\mathcal N (P_n)$, as well as the $(n+k)$-th Betti number of $\mathcal N (P_n)$ for $k \geq 0$. By way of Remark \ref{rem:8}, we get a lower bound for the dimension of the total cohomology of any finite simple graph.

We begin by formulating the asymptotic behavior of the total Betti numbers.

\begin{theorem}\label{thm:38}
As $n\to \infty$, the dimension $b(P_n,1)$ of the total cohomology $H^\bullet(\mathcal N(P_n))$ satisfies
\[
b(P_n,1)=2^{n-1}F_{n+2}\sim \frac{5+3\sqrt{5}}{20}(1+\sqrt{5})^n,
\]
where $F_{n}$ is the $n$-th Fibonacci number with $F_0=0$ and $F_1=1$.
\end{theorem}

\begin{proof}
Substituting $x=1$ in \eqref{eq:55}, we obtain
\[
b(P_n,1) = \frac{5+3\sqrt{5}}{20}(1+\sqrt{5})^n + \frac{5-3\sqrt{5}}{20}(1-\sqrt{5})^n,
\]
from which the result follows.
\end{proof}

\begin{rem} By work of the first author, Butler, Gardiner, Grandini, Lichtenwalner, and Pan \cite{ABGGLP}, the growth of the dimension of the total cohomology of paths is between that of totally disconnected graphs, $b(nK_1,1)=2^n$, and that of stars, $b(S_n,1)\sim \frac{4^{n+1}}{\sqrt{\pi n}}$. We can also compare our result to Tirao's work from \cite{Tirao00}, which gives the lower bound
\[b(P_n,1)\geq 2^{\dim(\mathfrak z (\mathcal L(P_n)))+\lfloor \frac{\dim(\mathfrak z (\mathcal L(P_n))^\perp)+1}{2}\rfloor}=2^{n-1+\lfloor\frac{n+1}{2}\rfloor},\]
where $\mathfrak z (\mathcal L(P_n))$ is the center of $\mathcal L(P_n)$ and $\mathfrak z (\mathcal L(P_n))^{\perp}$ is any complement of $\mathfrak z (\mathcal L(P_n))$. The gap between Tirao's lower bound and our formula is thus exponential because
\[
\lim_{n\to \infty} \sqrt[n]{\frac{b(P_n,1) }{2^{n-1+\lfloor\frac{n+1}{2}\rfloor}}}= \frac{1+\sqrt{5}}{2^{3/2}}\approx 1.14>1\,.
\]
\end{rem}

Theorem \ref{thm:38} gives a lower bound for the cohomology of Dani--Mainkar nilmanifolds.

\begin{cor}
For any finite simple graph $G$,
\[
\dim(H^\bullet(\mathcal N(G)))\ge 2^{{\rm diam}(G)}F_{{\rm diam}(G)+3}.
\]
\end{cor}

\begin{proof}
The claim follows from Theorem \ref{thm:38} combined with Remark \ref{rem:8}.
\end{proof}

The unimodality result in Proposition \ref{prop:38} implies that $b_n(P_n)$ is the largest Betti number of $\mathcal N(P_n)$. We conclude this paper with a description of the asymptotic behavior of a range of Betti numbers that includes $b_n(P_n)$.

\begin{theorem}\label{rangebetti}
Fix $k\geq 0$. As $n\to \infty$, 
\[
b_{n+k}(P_{n})\sim \sqrt{\frac{5+3\sqrt{5}}{40\pi n}}(1+\sqrt{5})^n.
\]
\end{theorem}

As our proof illustrates, the asymptotics in Theorem \ref{rangebetti} are independent of $k$.

\begin{proof}
By Poincar\'e duality, Cauchy's residue formula, and \eqref{eq:55},
\begin{align*}
b_{n+k}(P_n)&=b_{n-k-1}(P_n)\\
&=\frac{1}{2\pi \sqrt{-1}} \int_{C_1(0)} \frac{b_n(x)}{x^{n-k}} dx \\
&=\frac{1}{2\pi \sqrt{-1}} \int_{C_1(0)} \frac{A(x)}{x^{k}}e^{n\log(\psi_+(x)/x)} dx + \frac{1}{2\pi \sqrt{-1}} \int_{C_1(0)} \frac{B(x)}{x^k}e^{n\log(\psi_-(x)/x)} dx,
\end{align*}
where $C_1(0)$ denotes the unit circle in the complex plane. The critical points of $\log(\psi_\pm(x)/x)$ occur when $x\psi_{\pm}'(x)=\psi_\pm(x)$. Differentiating \eqref{eq:55bis}, we obtain
\[
\frac{\psi_\pm(x)}{x}=\psi_{\pm}'(x)=\frac{2x\psi_\pm+(1+x)(1+3x)}{2\psi_{\pm}(x)-1-x^2}
\]
whenever $x\neq 0$ is a critical point. By \eqref{eq:55bis}, we also get that $(1-x^2)(\psi_{\pm}(x)+x)=0$. Since $\psi_{\pm}(x)+x=0$ forces $x=0$, the only critical points on $C_{1}(0)$ are at $x=\pm 1$. Moreover, it follows from \eqref{eq:54} that, on $C_1(0)$,
\[
|\psi_{\pm}(x)|\le \frac{2+\sqrt{4+|1+x|^4}}{2}\le 1+\sqrt{5},
\]
with equality requiring $|1+x|=2$, that is, $x=1$. Then, the integrals involving $\psi_{\pm}$ contribute to at most $O((1+\sqrt{5}-\varepsilon)^n)$ on any closed circular segment of $C_1(0)$ that does not contain $x=1$. Since $\psi_-(1)=1-\sqrt{5}$, the integral involving $\psi_-$ contributes a total of $O((1+\sqrt{5}-\varepsilon)^n)$. It then follows from Laplace's method that
\[
b_{n+k}(P_n)\sim \frac{A(1)}{\sqrt{2\pi n C}} (1+\sqrt{5})^n,
\]
where $C$ is the result of evaluating the second derivative of $\log(\psi_+(x)/x)$ at $x=1$. Substituting \eqref{eq:54} yields
\[
C=A(1)=\frac{5+3\sqrt{5}}{20},
\]
from which the result follows.
\end{proof}

\bibliographystyle{abbrv}
\bibliography{Poincare110826}

\end{document}